\documentclass[12pt]{amsart}
\usepackage{amssymb}

\def\CC{{\mathbb C}}

\begin{document}

\title[Non-denseness of hyperbolicity]
{Non-denseness of hyperbolicity for linear isomorphisms in Banach spaces}
\author{Jos\'e F. Alves}
\address{Jos\'e Ferreira Alves, Centro de Matem\'atica da Universidade do Porto,
Rua do Campo Alegre 687, 4169-007 Porto, Portugal}
\email{jfalves@fc.up.pt}
\author{Maurizio Monge}
\address{Maurizio Monge,
Instituto de Matem\'atica da UFRJ,
Av. Athos da Silveira Ramos 149, Centro de Tecnologia,
Bloco C Cidade Universit\'aria,
Ilha do Fund\~ao, Caixa Postal 68530 21941-909, Rio de Janeiro, RJ, Brasil}
\email{maurizio.monge@im.ufrj.br}

\date{\today}

\thanks{JFA was partially supported by Funda\c c\~ao Calouste Gulbenkian, by ICTP, and by CMUP
  (UID/MAT/00144/2013), which is funded by FCT with national and European structural funds through the programs FEDER, under the partnership agreement PT2020.}
\thanks{MM was partially supported by ICTP and by EU Marie-Curie IRSES Brazilian-European partnership in Dynamical Systems (FP7-PEOPLE-2012-IRSES 318999 BREUDS)}



\begin{abstract}
We present an infinite dimensional Banach space in which  the set of hyperbolic linear isomorphisms in that space is not dense (in the norm topology) in the set of linear isomorphisms.
\end{abstract}

\maketitle

\section{Introduction}
The goal of  these note is to show that, contrarily to the finite dimensional case,  hyperbolicity is not necessarily dense in the space of linear isomorphisms of an infinite dimensional  Banach space. 

Consider $E$ a vector space endowed with a norm $\|\quad\|$. We say that $E$ is a  Banach space if the space $E$ with the metric induced by   $\|\quad\|$ is a complete metric space. We  say that a linear operator $T:E\to E$ is  bounded if 
 \begin{equation}\label{eq.norma}
\|T\|\equiv \sup_{x\neq 0}\frac{\|Tx\|}{\|x\|}=\sup_{\|x\|\le 1}{\|Tx\|}<\infty,
\end{equation}
and  define 
  $$L(E)=\{A:E\to E\mid\text{ $T$ is linear and bounded}\}.$$
 One has that $L(E)$ is a vector space and  $\|\quad\|$ defined as in~\eqref{eq.norma} gives a norm in $L(E)$.   Moreover,  if $E$ is a Banach space, then $L(E)$   is a Banach space as well.   
   Given $E$ a Banach space 
we denote by $GL(E)$ the set of invertible elements  in $L(E)$. Notice that as $E$ is a Banach space,   by  the Open Mapping Theorem  we   have $A^{-1}\in GL(E)$ whenever $A\in GL(E)$.

 We define the spectrum of
 $T\in GL(E)$  as
 $$\sigma(T)=\{\lambda\in\CC: \text{ $\lambda I-T$ is not invertible}\}.$$
We say that $A\in GL(E)$ is   \emph{hyperbolic} if the spectrum of $A$ is disjoint from the unit circle in $ \CC$, and define
 $$H(E)=\{A\in GL(E):\text{$A$ is hyperbolic}\}.$$
The dynamics of a hyperbolic linear isomorphism $A:E\to E$ is well understood in general: there are $A$-invariant linear subspaces $E^s $ and $E^u$ of $E$ with
\begin{enumerate}
\item  $E=E^s\oplus E^u$;
\item $A^nx\to 0$ as $n\to+\infty$ for all $x\in E^s$;
\item $A^{-n}x\to 0$ as $n\to+\infty$ for all $x\in E^u$;
\end{enumerate}
see e.g. \cite[Corollary~4.31]{I80}. Moreover, $H(E)$ is a dense subset of $GL(E)$; see e.g. \cite[Theorem 4.20]{I80}. 
In a finite dimensional space $E$ we have that $H(E)$ is also a dense subset of $GL(E)$; see e.g. \cite[Corollary~4.31]{I80}. 
It is then  natural to ask whether $H(E)$ is dense or not in $GL(E)$ when $E$ an infinite dimensional Banach space. In the next section we give an example of a Banach space $E$ where $H(E)$ is not dense in $GL(E)$.

\section{The example}

Let $D$ be the open disk of radius 1/2  centered at the point $1$ in the complex plane $\CC$, and let $S$ denote the boundary of $D$. 
Consider $ E$  the space of analytic functions $f: D\to\CC$ which have a continuous extension to~$S$, endowed with the supremum norm. By Morera's Theorem, we have that $E$ is a Banach space. 
Consider  $A:E\to E$ the bounded linear isomorphism in $E$ defined for each $f\in E$~as
 $$Af(z)=zf(z), \quad\forall z\in D.$$
As $A$ is a multiplication by the identity map in the domain $D$, we have  $\sigma(A)=\overline D$. We are going to see that $1\in \sigma (B)$ for any isomorphism $B$ sufficiently close to $A$ in the natural  norm of $ L( E)$.  Of course, this implies that $H( E)$ is not dense in $L(E)$. 
 
Consider  $\varepsilon>0$ and $B\in GL(E)$ with $\|A-B\|<\varepsilon$. To prove that  $1\in\sigma(B)$ for $B$ close to $A$, it is enough to see that $I-B$ is not surjective for $\varepsilon$ sufficiently small. Noticing that the constant function equal to 1 belongs to $E$, we shall actually prove that the equation 
  $(I-B)f=1$ has no solution $f\in E$. Notice that
$$ (I-B)f=1\quad\Leftrightarrow \quad f-Af=1-(A-B)f.$$
 Assume  by contradiction  that there is an $f\in E$ for which 
 \begin{equation}\label{eq.imposi}
f-Af=1-(A-B)f.
\end{equation}
 Using the Maximum Modulus Principle we may write
\begin{eqnarray*}
\|f-Af\|&=&\sup_{z\in D}|f(z)-zf(z)|\\
&=&\sup_{z\in S}|f(z)-zf(z)|\\
&=&\sup_{z\in S}|(1-z)f(z)|\\
&=&\frac12\|f\|,
\end{eqnarray*}
which then yields
 $$\|(A-B)f\|\le\varepsilon\|f\|=2\varepsilon \|f-Af\|,$$
 and so
 $$ \|f-Af\|\ge \frac1{2\varepsilon}\|(A-B)f\|.$$
Using~\eqref{eq.imposi} we obtain
 $$\frac1{2\varepsilon}\|(A-B)f\|\le  \|f-Af\|\le   1+\|(A-B)f\|.$$
Hence, for $\varepsilon<1/2$ we have
  $$\|(A-B)f\|\le \left(\frac1{2\varepsilon}-1\right)^{-1}<1.$$
  This in particular gives that the function $1-(A-B)f$ has no zeros in $D$ by Rouch\'e Theorem
  applied to the functions $1$ and $(A-B)f$.
 Noticing that the equation $(f-Af)(z)=0$ has a solution in $z=1$, we have a contradiction with \eqref{eq.imposi}.


\begin{thebibliography}{1}

\bibitem{I80}
M.~C. Irwin.
\newblock {\em Smooth dynamical systems}, volume~94 of {\em Pure and Applied
  Mathematics}.
\newblock Academic Press Inc. [Harcourt Brace Jovanovich Publishers], New York,
  1980.

\end{thebibliography}

\end{document}